# ON THE WEAK-FRAGMENTABILITY INDEX OF SOME LIPSCHITZ-FREE SPACES

ESTELLE BASSET

ABSTRACT. We show the existence of Lipschitz-free spaces satisfying the Point of Continuity Property with arbitrarily high weak-fragmentability index. For this purpose, we use a generalized construction of the countably branching diamond graphs. A direct corollary is the existence of an uncountable family of pairwise non-isomorphic Lipschitz-free spaces over purely 1-unrectifiable metric spaces, but more importantly, the existence of a Lipschitz-free space over a separable metric space that is not isomorphic to any Lipschitz-free space over a compact. Another consequence is that to be Lipschitz-universal for countable complete metric spaces, a separable complete metric space cannot be purely 1-unrectifiable. Some results on compact reduction are also obtained.

## 1. INTRODUCTION

Given a metric space $M$ with a distinguished base point $0 \in M$, we will denote by $\mathrm{Lip}_0(M)$ the vector space of Lipschitz functions $f\colon M \to \mathbb{R}$ vanishing at the designated origin 0. We will endow $\mathrm{Lip}_0(M)$ with the norm given by the best Lipschitz constant

$$\|f\|_L = \sup\left\{\frac{f(x)-f(y)}{d(x,y)},\ x \neq y \in M\right\}$$

(which is not a norm in the vector space of Lipschitz functions, that is why we work with the space $\mathrm{Lip}_0(M)$ instead); equipped with this norm, $\mathrm{Lip}_0(M)$ turns out to be a Banach space. For $x \in M$, we let $\delta_M(x) \in \mathrm{Lip}_0(M)^*$ be the evaluation functional $\delta_M(x)\colon f \mapsto f(x)$, and we define the *Lipschitz-free space over* $M$, denoted by $\mathcal{F}(M)$, as the norm-closure of the linear span of $\{\delta_M(x),\ x \in M\}$ in $\mathrm{Lip}_0(M)^*$. It is readily seen that $\delta_M$ is an isometric embedding of $M$ into $\mathcal{F}(M)$.

Recall that a metric space is *purely* 1*-unrectifiable* if and only if it contains no bi-Lipschitz copies of compact, positive measure subsets of $\mathbb{R}$. This kind of metric space is of significant interest when dealing with Lipschitz-free spaces, considering the following theorem due to Aliaga, Gartland, Petitjean and Procházka in [2]:

**Theorem.** *Let $M$ be a metric space. The following are equivalent:*

*(i) The completion of $M$ is purely* 1*-unrectifiable.*
*(ii) $\mathcal{F}(M)$ has the Radon–Nikodým property (RNP).*
*(iii) $\mathcal{F}(M)$ has the Schur property.*
*(iv) $\mathcal{F}(M)$ has the Krein–Milman property.*
*(v) $\mathcal{F}(M)$ does not contain any isomorphic copy of $L_1$.*

---

We recall that a Banach space $X$ is said to have the *Point of Continuity Property* (PCP) if every non-empty bounded closed subset $F$ of $X$ has a weak point of continuity, that is, for every such subset $F$ of $X$, the identity map

$$id\colon (F,\, w) \to (F,\, \|.\|_X)$$

is continuous at some point $x \in F$ (where $w$ stands for the weak topology on $F$). In geometric terms (see for example Definition 2.25 in [23]), $X$ has the PCP if and only if every non-empty bounded subset of $X$ is *w-fragmentable*, *i.e.* has non-empty relatively $w$-open subsets of arbitrarily small diameter. Since a Banach space has the RNP if and only if every non-empty bounded subset of $X$ is *dentable*, that is, has open slices of arbitrarily small diameter, the RNP implies the PCP. As well, $L_1$ does not have the PCP, so the PCP is an intermediate property between the RNP and non-containment of $L_1$. Therefore, we have the following characterization (see Remark 4.7 in [2]):

**Remark 1.1.** A metric space $M$ has a purely 1-unrectifiable completion if and only if $\mathcal{F}(M)$ has the PCP.

In this paper, we will use an ordinal $\Phi(X)$, called the *weak-fragmentability index of* $X$, to testify about "how much" a Banach space $X$ has the PCP. Actually, $\Phi$ is the weak version of the Szlenk index, introduced by Szlenk in [25]. We will recall its definition and some basic properties in Section 2.2. If $M$ is a purely 1-unrectifiable separable complete metric space, the weak-fragmentability index of its Lipschitz-free space is strictly less than $\omega_1$ (first uncountable ordinal, see Proposition 2.5). But we show that there exist such metric spaces, and even countable complete metric spaces, whose Lipschitz-free space has arbitrarily large weak-fragmentability index:

**Theorem** (Theorem 3.7)**.** *For every $\alpha \in (0,\, \omega_1)$, there exists a countable complete metric space $D_\alpha$ such that $\Phi(\mathcal{F}(D_\alpha)) > \alpha$.*

The construction of the $D_\alpha$'s is given in Section 2.3, the proof of Theorem 3.7 is in Section 3.

We develop some consequences in Section 4, in particular the following one about universality:

**Corollary** (Corollary 4.1)**.** *Let $M$ be a separable complete metric space such that every countable complete metric space is Lipschitz-equivalent to a subspace of $M$. Then $M$ is not purely 1-unrectifiable.*

The $D_\alpha$ spaces will also provide a negative answer to a previously open question:

**Question 1.2.** Given a separable metric space $M$, does there exist a compact space $K$ such that $\mathcal{F}(M)$ is linearly isomorphic to $\mathcal{F}(K)$?

Indeed, this is not the case for the $D_\alpha$ spaces with $\alpha \in [\omega, \omega_1)$, see Corollary 4.3.

We say a word about compact reduction in Proposition 4.7. Finally, we raise a question related to the dentability index $D$ (see Section 5):

**Question** (Question 5.2)**.** Does there exist a function $\Psi\colon (0,\, \omega_1) \to (0,\, \omega_1)$ such that for every separable metric space $M$ and for any $\alpha \in (0,\, \omega_1)$, $\Phi(\mathcal{F}(M)) \leq \alpha$ implies $D(\mathcal{F}(M)) \leq \Psi(\alpha)$?

## 2. Background

In this section, we introduce the main objects used in this paper, set some notation and recall some elementary properties.

2.1. **Lipschitz-free spaces.** To simplify the notation we write $\delta$ instead of $\delta_M$ when there is no ambiguity. By a *molecule* we mean an element of $\mathcal{F}(M)$ of the form

$$m_{x,y} = \frac{\delta(x) - \delta(y)}{d(x,y)}$$

for $x \neq y \in M$. Notice that molecules are of norm one. These elements will play a key role in the sequel.

A fundamental property of $\mathcal{F}(M)$ is the following "universal extension property": any Banach-space valued Lipschitz function $f\colon M \to X$ vanishing at 0 can be uniquely extended (identifying $M$ with $\delta(M)$) to a continuous linear map $\widehat{f}\colon \mathcal{F}(M) \to X$ whose operator norm is equal to $\|f\|_L$. An easy consequence of this property, picking $X = \mathbb{R}$, is that $\mathcal{F}(M)^*$ is linearly isometric to $\mathrm{Lip}_0(M)$. Another useful observation is that whenever $N$ is a subset of $M$, then $\mathcal{F}(N)$ is linearly isometric to a subspace of $\mathcal{F}(M)$. More precisely:

**Proposition 2.1.** *If $N$ is a subset of $M$ containing $0$, the map $\iota_N$ defined by*

$$\iota_N\left(\sum_{i=1}^{n} a_i \delta_N(x_i)\right) = \sum_{i=1}^{n} a_i \delta_M(x_i),\ \ x_1,\dots,x_n \in N,\ \ a_1,\dots,a_n \in \mathbb{R}$$

*can be extended to a linear isometry from $\mathcal{F}(N)$ to $\mathcal{F}_N(M)$, the closed linear span of $\delta_M(N)$ in $\mathcal{F}(M)$.*

From now on we will use freely this identification.

Two metric spaces $M$, $N$ are said to be *Lipschitz-equivalent* if there exists a bijection $f\colon M \to N$ such that $f$ and $f^{-1}$ are Lipschitz maps. Using the previous identification, the universal extension property enables us to prove that if $N$ and $M$ are Lipschitz-equivalent, then $\mathcal{F}(N)$ and $\mathcal{F}(M)$ are linearly isomorphic. For a quick proof of the universal extension property and some other basic facts about Lipschitz-free spaces, we refer the reader to [8].

2.2. **The weak-fragmentability index.** Now we give the definition of the weak-fragmentability index and the matching derivation, and review its basic properties. In the sequel, all Banach spaces we consider are over the real field.

Let $X$ be a Banach space and $K$ be a $w$-closed bounded subset of $X$. For every $\varepsilon > 0$, we define the *derived set of $K$*

$$\sigma_\varepsilon(K) \coloneqq K \setminus \bigcup \{V \subset X \ w-\mathit{open} : \mathrm{diam}(V \cap K) < \varepsilon\}.$$

In other words, we remove from $K$ all points admitting relative weak-neighborhood of diameter less than $\varepsilon$. Then, given an ordinal $\alpha$, we define inductively $\sigma_\varepsilon^\alpha(K)$ by setting $\sigma_\varepsilon^0(K) = K$, $\sigma_\varepsilon^{\alpha+1}(K) = \sigma_\varepsilon(\sigma_\varepsilon^\alpha(K))$ and $\sigma_\varepsilon^\alpha(K) = \bigcap_{\beta<\alpha} \sigma_\varepsilon^\beta(K)$ if $\alpha$ is a limit ordinal.

We denote by $B_X$ the closed unit ball of $X$. We then define $\Phi(X,\varepsilon)$ as the smallest ordinal $\alpha$ such that $\sigma_\varepsilon^\alpha(B_X) = \varnothing$, when such an ordinal exists. Otherwise, we write $\Phi(X,\varepsilon)=\infty$. If $\Phi(X,\varepsilon)$ is defined for all $\varepsilon>0$, then we define the *weak-fragmentability index of* $X$ as the ordinal

$$\Phi(X) \coloneqq \sup_{\varepsilon>0}\Phi(X,\varepsilon),$$

and we write $\Phi(X)<\infty$ to signify that $\Phi(X)$ is well defined. If $\Phi(X,\varepsilon)=\infty$ for some $\varepsilon>0$, we write $\Phi(X)=\infty$.

First, let us give two elementary facts: if $X$ is linearly isomorphic to another Banach space $Y$, then $\Phi(X)=\Phi(Y)$, and if $F$ is a closed subspace of $X$, then $\Phi(F)\le\Phi(X)$. Now we show two useful properties. Actually, Proposition 2.2 and Fact 2.4 below are proved in [20] (Proposition 2) in the weak$^*$ case (for the Szlenk index). See [5] (Proposition II.2) for Proposition 2.5 below. Nevertheless, we reproduce the proofs here for the reader's convenience.

**Proposition 2.2.** *Let $X$ be a Banach space. If $\Phi(X)<\infty$, then there exists an ordinal $\alpha$ such that $\Phi(X)=\omega^\alpha$ (where $\omega$ denotes the first infinite ordinal).*

In order to show this, we will need two facts:

**Fact 2.3.** *For every ordinal $\alpha$, for every $\varepsilon>0$, we have the inclusion:*

$$\frac{1}{2}\sigma_\varepsilon^\alpha(B_X)+\frac{1}{2}B_X\subset\sigma_{\varepsilon/2}^\alpha(B_X).$$

This fact follows easily with a transfinite induction on $\alpha$:

*Proof.* The above inclusion is clear if $\alpha$ equals 0. Assume it is true for every $\beta<\alpha$. If $\alpha$ is a limit ordinal, we get the conclusion by taking the intersection on $\beta<\alpha$. If $\alpha=\beta+1$ is a successor ordinal, let $z=\frac{1}{2}x+\frac{1}{2}y$ with $x\in\sigma_\varepsilon^\alpha(B_X)$ and $y\in B_X$. By induction hypothesis, $z\in\sigma_{\varepsilon/2}^\beta(B_X)$. Let $V$ be a relative $w$-open neighborhood of $z$ in $\sigma_{\varepsilon/2}^\beta(B_X)$. We must show that $\operatorname{diam}(V)\ge\frac{\varepsilon}{2}$. Without loss of generality, we may assume that $V$ is of the form

$$V=\{z'\in\sigma_{\varepsilon/2}^\beta(B_X):\forall r\in\{1,\dots,n\},\ \left|\langle f_r,\ z-z'\rangle\right|\le\eta\}$$

with $\eta>0$ and $f_1,\dots,f_n\in X^*$. By induction hypothesis, $V$ contains the set $\frac{1}{2}W+\frac{1}{2}y$ where

$$W=\{x'\in\sigma_\varepsilon^\beta(B_X):\forall r\in\{1,\dots,n\},\ \left|\langle f_r,\ x-x'\rangle\right|\le\eta\}$$

is a relative $w$-open neighborhood of $x$ in $\sigma_\varepsilon^\beta(B_X)$. Thus $\operatorname{diam}(W)\ge\varepsilon$, and therefore

$$\operatorname{diam}(V)\ge\operatorname{diam}(\frac{1}{2}W+\frac{1}{2}y)\ge\frac{\varepsilon}{2}.$$

□

**Fact 2.4.** *For every ordinal $\alpha$, we have*

$$\Phi(X)>\omega^\alpha\Longrightarrow\Phi(X)\ge\omega^{\alpha+1}.$$

*Proof.* Assume $\Phi(X)=\sup_{\varepsilon>0}\Phi(X,\varepsilon)>\omega^\alpha$. Then there exists $\varepsilon>0$ such that $\Phi(X,2\varepsilon)>\omega^\alpha$ *i.e.* $\sigma_{2\varepsilon}^{\omega^\alpha}(B_X)\ne\varnothing$: let $x$ be an element of this set. Since $-x\in B_X$, Fact 2.3 implies that $0=\frac{1}{2}x+\frac{1}{2}(-x)\in\sigma_\varepsilon^{\omega^\alpha}(B_X)$. Using Fact 2.3 again, we get that $\frac{1}{2}B_X\subset\sigma_{\varepsilon/2}^{\omega^\alpha}(B_X)$, and consequently $\sigma_{\varepsilon/2}^{\omega^\alpha}(\frac{1}{2}B_X)\subset\sigma_{\varepsilon/2}^{\omega^\alpha.2}(B_X)$. Recalling that $0\in\sigma_\varepsilon^{\omega^\alpha}(B_X)$, we then have $0\in\sigma_{\varepsilon/2}^{\omega^\alpha}(\frac{1}{2}B_X)\subset\sigma_{\varepsilon/2}^{\omega^\alpha.2}(B_X)$. More generally, we can show with an induction that for every $n\in\mathbb{N}$, $0\in\sigma_{\varepsilon/2^n}^{\omega^\alpha.2^n}(B_X)$, so the latter set is non-empty. Hence, $\Phi(X) \geq \Phi(X, \frac{\varepsilon}{2^n}) > \omega^\alpha . 2^n$ for all $n$, so $\Phi(X) \geq \omega^\alpha . \omega = \omega^{\alpha+1}$. □

We are now able to prove the desired proposition:

*Proof of Proposition 2.2.* Let $\alpha$ be the infimum of all ordinals $\gamma$ such that $\Phi(X) \leq \omega^\gamma$, well defined since $\Phi(X) < \infty$. If $\alpha$ is a limit ordinal, $\omega^\alpha = \sup_{\beta<\alpha} \omega^\beta \leq \Phi(X)$ since $\omega^\beta < \Phi(X)$ for every $\beta < \alpha$, and then $\omega^\alpha = \Phi(X)$ by definition of $\alpha$. If $\alpha = \beta + 1$, then $\Phi(X) > \omega^\beta$ and Fact 2.4 leads to $\Phi(X) \geq \omega^{\beta+1} = \omega^\alpha$. As $\Phi(X) \leq \omega^\alpha$, we finally have $\Phi(X) = \omega^\alpha$. □

The PCP is related to the weak-fragmentability index in the following well-known way:

**Proposition 2.5.** *Let $X$ be a separable Banach space. Then $X$ has the PCP if and only if $\Phi(X) < \omega_1$.*

*Proof.* Assume $X$ has the PCP, and let $\varepsilon > 0$. For every ordinal $\alpha < \Phi(X, \varepsilon)$, the set $\sigma_\varepsilon^\alpha(B_X)$ is non-empty. But $\sigma_\varepsilon^\alpha(B_X)$ is a bounded closed subset of $X$ which has the PCP, so it admits a point of continuity: $\sigma_\varepsilon^\alpha(B_X) \setminus \sigma_\varepsilon^{\alpha+1}(B_X) \neq \varnothing$. Thus, there exists an open set $\mathcal{O}$ such that $\mathcal{O} \cap \sigma_\varepsilon^\alpha(B_X) \neq \varnothing$ and $\mathcal{O} \cap \sigma_\varepsilon^{\alpha+1}(B_X) = \varnothing$. Since $X$ is separable, it has a countable open base $(\mathcal{O}_n)_n$ and $\mathcal{O}$ can be written as the union of some $\mathcal{O}_n$'s. Consequently, there exists $n_\alpha \in \mathbb{N}$ such that $\mathcal{O}_{n_\alpha} \cap \sigma_\varepsilon^\alpha(B_X) \neq \varnothing$ and $\mathcal{O}_{n_\alpha} \cap \sigma_\varepsilon^{\alpha+1}(B_X) = \varnothing$. Next, define $f$ a function mapping each $\alpha < \omega_1$ to such an $n_\alpha$. By definition of $n_\alpha$, $f$ is a one-to-one mapping from $(0, \Phi(X, \varepsilon))$ to $\mathbb{N}$, which implies $\Phi(X, \varepsilon) < \omega_1$.

Due to the monotonicity of $\Phi(X, \varepsilon)$ with respect to $\varepsilon$, we have the equality $\sup_{\varepsilon>0} \Phi(X, \varepsilon) = \sup_{n\in\mathbb{N}} \Phi(X, \frac{1}{n})$. Hence, $\Phi(X)$ is a countable supremum of countable ordinals, so is countable itself.

For the converse, if $X$ does not have the PCP, there are $B \neq \varnothing$ a subset of $B_X$ and $\varepsilon > 0$ such that $\sigma_\varepsilon(B) = B$. Therefore, $\sigma_\varepsilon^\alpha(B) = B \neq \varnothing$ for all $\alpha$, so $\Phi(X) = \infty$. □

**Remark 2.6.** The weak-fragmentability index has already been used under the name of "weak Szlenk index". The first occurence we could find was in [19]. This index is less known than the classical Szlenk index because, due in particular to the lack of weak-compactness, we are still lacking results as important as those related to the Szlenk index, for example in renorming theory or non-linear geometry of Banach spaces. Let us mention that it is still an important open question to know whether a Banach space $X$ satisfies $\Phi(X) \leq \omega$ if and only if it is AUC renormable (see Section 4.3 for the definition of AUC renormability). We have chosen to call it "weak fragmentability index", as we felt that it is more related to the fundamental works of Jayne and Rogers [15] and Namioka [21] than to the work of Szlenk. Indeed, in [15] and [21], a metric space $(M, d)$ is said to be *fragmentable by a topology* $\tau$ on $M$ if for every $\varepsilon > 0$ and every non empty $\tau$-closed subset $F$ of $M$, there exists $U \in \tau$ such that $U \cap F \neq \emptyset$ and the diameter (for the metric $d$) of $U \cap F$ is less than $\varepsilon$.

2.3. **The countably branching diamond graphs.** We will define a transfinite sequence of metric spaces $(D_\alpha, d_\alpha)$ together with distinguished distinct points $t_\alpha$, $b_\alpha$, $\ell_\alpha \in D_\alpha$ (informally speaking: the top, the bottom and the left-most points of $D_\alpha$) such that $d_\alpha(t_\alpha, b_\alpha) = 2$, indexed by $\alpha \in (0, \omega_1)$.

In order to best describe the metric of $D_\alpha$, we find it convenient to consider $(D_\alpha, d_\alpha)$ as a weighted graph. In particular, we describe only the distances of points connected by an edge, while all the other distances are computed as the shortest path distance in the corresponding weighted graph.

We start by defining $D_1 = \{t_1, b_1\} \cup \{x_1^n : n \in \mathbb{N}\}$. We declare $\ell_1 := x_1$ and the only edges will be the couples of the form $\{t_1, x_1^n\}$ and $\{b_1, x_1^n\}$ for some $n \in \mathbb{N}$. All the edges have weight 1. Notice that $d_1(t_1, b_1) = 2$.

Now assume that $(D_\beta, d_\beta)$ together with $t_\beta$, $b_\beta$, $\ell_\beta \in D_\beta$ have been constructed and also the graph structure on $D_\beta$ is known, and let $\alpha = \beta + 1$. In order to obtain $D_\alpha$, we will replace each edge of $D_1$ by a copy of $D_\beta$ scaled-down by the factor 2. More precisely, let $D_\alpha^{(j,+)}$, resp. $D_\alpha^{(j,-)}$, $j \in \mathbb{N}$, be pairwise disjoint isometric copies of $(D_\beta, \frac{d_\beta}{2})$ (in particular, in the corresponding weighted graph we also divide the weights by 2) and let $t_\alpha^{(j,+)}$, $b_\alpha^{(j,+)}$, resp. $t_\alpha^{(j,-)}$, $b_\alpha^{(j,-)}$, be the elements of $D_\alpha^{(j,+)}$, resp. $D_\alpha^{(j,-)}$, corresponding respectively to $t_\beta$, $b_\beta$ through this isometry. We define

$$D_\alpha := \dot{\bigcup_{j\in\mathbb{N}}} D_\alpha^{(j,+)} \,\dot\cup\, D_\alpha^{(j,-)},$$

where the dot symbol over the union is only here to signify that, in the union $\bigcup_{j\in\mathbb{N}} D_\alpha^{(j,+)} \cup D_\alpha^{(j,-)}$, we identify or "glue together" some points:

- $b_\alpha^{(j,+)}$ with $t_\alpha^{(j,-)}$, for all $j \in \mathbb{N}$ (we glue $D_\alpha^{(j,+)}$ with $D_\alpha^{(j,-)}$ by one of their poles);
- the $t_\alpha^{(j,+)}$ points all together, for $j \in \mathbb{N}$ (we glue the sets $D_\alpha^{(j,+)} \,\dot\cup\, D_\alpha^{(j,-)}$ by their top poles) and we denote $t_\alpha$ the resulting point;
- the $b_\alpha^{(j,-)}$ points all together, for $j \in \mathbb{N}$ (we glue the sets $D_\alpha^{(j,+)} \,\dot\cup\, D_\alpha^{(j,-)}$ by their bottom poles) and we denote $b_\alpha$ the resulting point.

The set of edges in the graph $D_\alpha$ is the union of the edges of all graphs $D_\alpha^{(j,-)}$ and $D_\alpha^{(j,+)}$, together with their weights. Also, we declare $\ell_\alpha := t_\alpha^{1,-} = b_\alpha^{1,+}$. Observe that $d_\alpha(t_\alpha, b_\alpha) = 2$.

Finally, if $\alpha < \omega_1$ is a limit ordinal and if for each $\beta < \alpha$, $D_\beta$ and the other objects have been defined as wanted, we glue the $D_\beta$ graphs together along the sets $\{t_\beta\}_{\beta<\alpha}$, $\{b_\beta\}_{\beta<\alpha}$. More precisely, we set

$$D_\alpha := \dot{\bigcup_{\beta<\alpha}} \{\beta\} \times D_\beta$$

where again, the dot indicates that we consider some points to be the same:

- all top poles $(\beta, t_\beta)$ together, and we denote the resulting point $t_\alpha$;
- all bottom poles $(\beta, b_\beta)$ together, and we denote the resulting point $b_\alpha$.

This construction is valid since by induction hypothesis, $d_\beta(t_\beta, b_\beta) = 2$ for each $\beta < \alpha$. So $d_\alpha(t_\alpha, b_\alpha) = 2$. Also, we declare $\ell_\alpha := \ell_1$. In $D_\alpha$, there is an edge between two elements $(\beta, x)$, $(\gamma, y)$ if and only if $\beta = \gamma$ and there is an edge between $x$ and $y$ in $D_\beta$. In this case, the weight on this edge is the same as the weight on the corresponding $D_\beta$'s edge. This ends the induction.

FIGURE 1. $D_1$ and $D_\alpha$ for $\alpha$ a successor

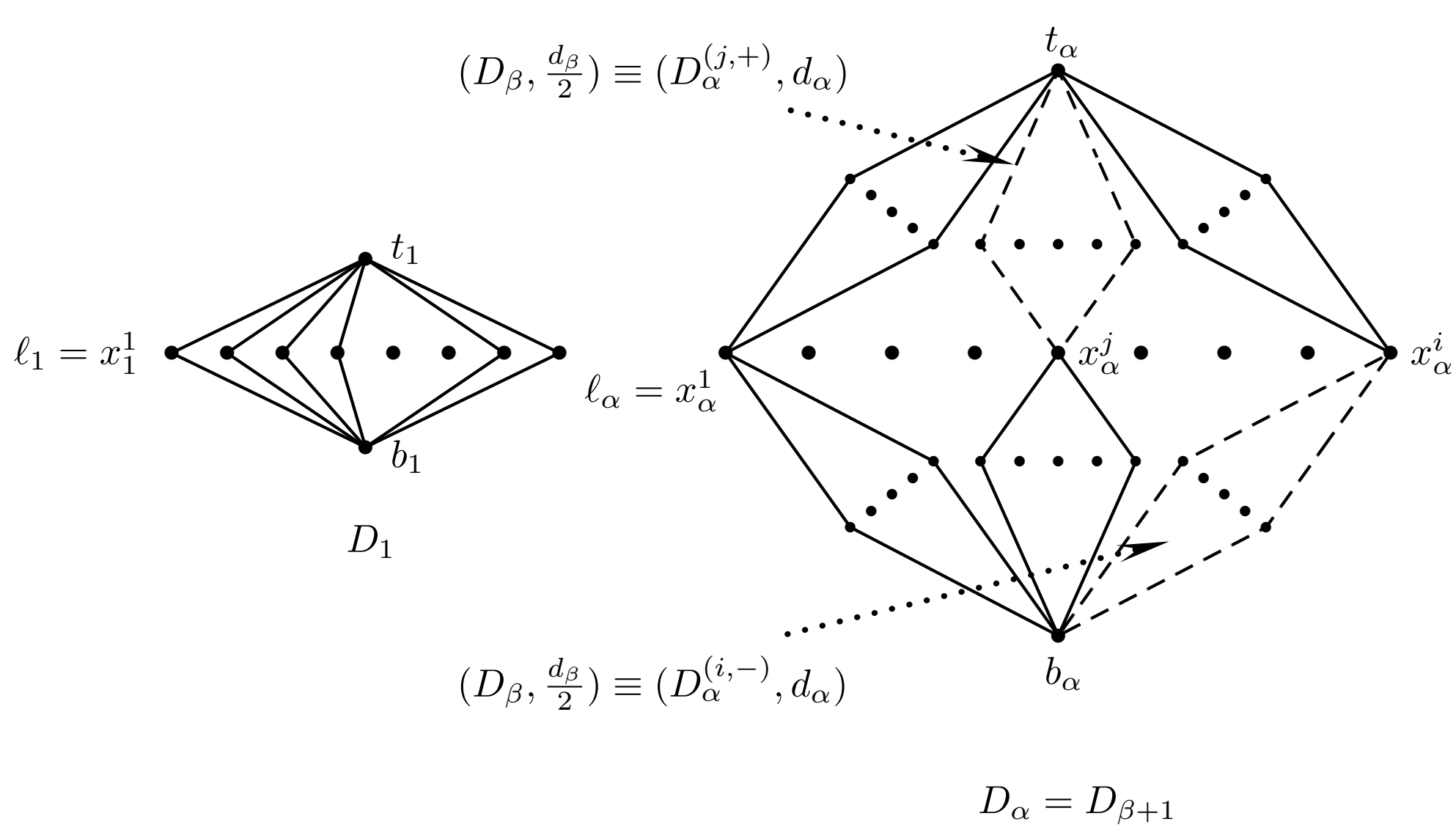

This construction is a generalization of the classical sequence $(D_k)_{k\in\mathbb{N}}$ of diamond graphs. For further details on these graphs, we refer the reader to [4].

In the subsequent proofs, we will have to address the "middle points" of $D_\alpha$ when $\alpha$ is a successor ordinal. So we denote by $(x_\alpha^n)_{n\in\mathbb{N}}$ the points of $D_\alpha$ corresponding to the $x_1^n$'s in $D_1$ (see Figure 1): $x_\alpha^j = t_\alpha^{(j,-)} = b_\alpha^{(j,+)}$.

For every $\alpha \in (0, \omega_1)$, we distinguish $\ell_\alpha$ as the base point of $D_\alpha$, and we will be interested in the Lipschitz-free space on $D_\alpha$. First, let us notice that $\mathcal{F}(D_\alpha)$ has the PCP (see Remark 1.1):

**Proposition 2.7.** *For every $\alpha \in (0, \omega_1)$, $D_\alpha$ is a countable complete metric space, and thus a purely 1-unrectifiable metric space.*

*Proof.* The countability of the $D_\alpha$'s spaces is clear by transfinite induction. We prove the completeness of the $D_\alpha$'s with a transfinite induction too. $D_1$ is complete since it is uniformly discrete. Assume now that $D_\beta$ is complete for every $\beta < \alpha$. Let $(x_n)_n \subset D_\alpha$ be a Cauchy sequence. We consider several cases, each time using the fact that a Cauchy sequence which has a convergent subsequence is convergent itself.

Let us first assume that $\alpha$ is a limit ordinal.

Case 1: the terms of the sequence belong to a finite number of $D_\beta$'s. Then there exist $\beta < \alpha$ and a subsequence of $(x_n)_n$ included in $D_\beta$. But $D_\beta$ is complete by induction hypothesis, so $(x_n)_n$ is convergent.

Case 2: the terms of the sequence belong to an infinite number of $D_\beta$'s and for every $\varepsilon > 0$, there exists $n \in \mathbb{N}$ such that $\mathrm{dist}(x_n, \{t_\alpha, b_\alpha\}) < \varepsilon$. Then there exist a subsequence of $(x_n)_n$ which converges towards one of the poles, so $(x_n)_n$ is convergent.

Case 3: the terms of the sequence belong to an infinite number of $D_\beta$'s and there exists $\varepsilon > 0$ such that for every $n \in \mathbb{N}$, $\operatorname{dist}(x_n, \{t_\alpha, b_\alpha\}) \geq \varepsilon$. But to trace a path between two elements belonging to different $D_\beta$'s, we have to pass through one of the poles, so there exists a subsequence of $(x_n)_n$ which is $2\varepsilon$-separated. This is in contradiction with the fact that $(x_n)_n$ is a Cauchy sequence.

If $\alpha$ is a successor ordinal, we reproduce the same proof by exhaustion with the spaces $D_\alpha^{(i,+)} \cup D_\alpha^{(i,-)}$, $i \in \mathbb{N}$ (which are complete by induction hypothesis) playing the role of the spaces $D_\beta$, $\beta < \alpha$, which ends the proof. □

## 3. Minoration of $\Phi(\mathcal{F}(D_\alpha), 1)$

In this section, we will show that there exists an $\varepsilon > 0$ such that for every $\alpha < \omega_1$, the set $\sigma_\varepsilon^\alpha(B_{\mathcal{F}(D_\alpha)})$ is non-empty.

### 3.1. First derived set of $B_{\mathcal{F}(D_\alpha)}$.

We start by noticing that the molecule associated with the two poles of $\mathcal{F}(D_\alpha)$ is in the first derived set of $B_{\mathcal{F}(D_\alpha)}$ for $\varepsilon = 1$.

**Proposition 3.1.** *For every $\alpha \in (0, \omega_1)$, $m_{t_\alpha, b_\alpha} \in \sigma_1(B_{\mathcal{F}(D_\alpha)})$.*

*Proof.* Let first assume that $\alpha$ is equal to 1 or a successor ordinal. Let $V$ be a $w$-open neighborhood of $m_{t_\alpha,b_\alpha}$ in $\mathcal{F}(D_\alpha)$. We must show that $\operatorname{diam}(V \cap B_{\mathcal{F}(D_\alpha)}) \geq 1$. Without loss of generality, we may assume that $V$ is of the form

$$V = \{\mu \in \mathcal{F}(D_\alpha) : \forall r \in \{1, \dots, n\},\ |\langle f_r,\ \mu - m_{t_\alpha,b_\alpha}\rangle| \leq \varepsilon\}$$

with $\varepsilon > 0$ and $f_1, \dots, f_n \in \operatorname{Lip}_0(D_\alpha)$. Since $d_\alpha(t_\alpha, b_\alpha) = 2d_\alpha(t_\alpha, x_\alpha^i) = 2d_\alpha(x_\alpha^i, b_\alpha)$ for all $i \in \mathbb{N}$, we can write:

$$m_{t_\alpha,b_\alpha} = \frac{1}{2}(\underbrace{m_{t_\alpha,x_\alpha^i}}_{:=\mu_i} + \underbrace{m_{x_\alpha^i,b_\alpha}}_{:=\nu_i}).$$

Let us pass to some subsequences $(\mu_{i_j})_{j\in\mathbb{N}}$ and $(\nu_{i_j})_{j\in\mathbb{N}}$ such that for every $r \in \{1, \dots, n\}$, the sequences $(f_r(\mu_{i_j}))_j$ and $(f_r(\nu_{i_j}))_j$ are convergent in $\mathbb{R}$. Thus, we have, for every $r \in \{1, \dots, n\}$:

$$\begin{aligned}\langle f_r,\ m_{t_\alpha,b_\alpha}\rangle &= \left\langle f_r,\ \frac{\mu_{i_j} + \nu_{i_j}}{2}\right\rangle \text{ for all } j \in \mathbb{N} \\ &= \lim_{j\to+\infty} \left\langle f_r,\ \frac{\mu_{i_j} + \nu_{i_j}}{2}\right\rangle \\ &= \lim_{j\to+\infty} \left\langle f_r,\ \frac{\mu_{i_j+1} + \nu_{i_j}}{2}\right\rangle.\end{aligned}$$

Therefore, there is $j \in \mathbb{N}$ large enough such that $\gamma_V := \frac{\mu_{i_j+1}+\nu_{i_j}}{2} \in V$. Given such a $j$, we have that

$$\begin{aligned}\|\gamma_V - m_{t_\alpha,b_\alpha}\|_{\mathcal{F}(D_\alpha)} &= \left\|\frac{\mu_{i_j+1} - \mu_{i_j}}{2}\right\| \\ &= \frac{1}{2}\left\|\frac{\delta(t_\alpha) - \delta(x_\alpha^{i_j+1})}{d_\alpha(t_\alpha, x_\alpha^{i_j+1})} - \frac{\delta(t_\alpha) - \delta(x_\alpha^{i_j})}{d_\alpha(t_\alpha, x_\alpha^{i_j})}\right\| \\ &= \left\|\frac{\delta(x_\alpha^{i_j}) - \delta(x_\alpha^{i_j+1})}{2}\right\| = \left\|m_{x_\alpha^{i_j}, x_\alpha^{i_j+1}}\right\| = 1,\end{aligned}$$

so $\operatorname{diam}(V \cap B_{\mathcal{F}(D_\alpha)}) \geq 1$.

Assume now that $\alpha$ is a limit ordinal. Then $m_{t_\alpha,b_\alpha}$ identifies with $m_{t_\beta,b_\beta} \in \mathcal{F}(D_\beta)$ for some successor ordinal $\beta < \alpha$. By the previous arguments, $m_{t_\beta,b_\beta} \in \sigma_1(B_{\mathcal{F}(D_\beta)})$, and $\sigma_1(B_{\mathcal{F}(D_\beta)}) \subsetneqq \sigma_1(B_{\mathcal{F}(D_\alpha)})$, so $m_{t_\alpha,b_\alpha} \in \sigma_1(B_{\mathcal{F}(D_\alpha)})$. □

**Remark 3.2.** Notice that in the process, we proved that for every successor ordinal $\alpha$ and for every $V$ $w$-open neighborhood of $m_{t_\alpha,b_\alpha}$ in $\mathcal{F}(D_\alpha)$, there exist $j > i$ in $\mathbb{N} \setminus \{1\}$ such that $\gamma_V \coloneqq \frac{1}{2}(m_{t_\alpha,x_\alpha^j} + m_{x_\alpha^i,b_\alpha}) \in V$.

3.2. **Stability by taking special convex combinations.** We show that in some configurations, if two elements belong to a derived set of $B_{\mathcal{F}(D_\alpha)}$, then their average also belongs to it.

Recall that $\ell_\alpha$ is the base point of $D_\alpha$. For $\alpha$ a successor ordinal, we use the pieces of notation $D_\alpha^{(j,+)}$ and $D_\alpha^{(i,-)}$ introduced in Section 2.3. Each of these spaces is considered along with its origin $\ell_\alpha^{(j,+)}$ and $\ell_\alpha^{(i,-)}$ corresponding to the origin $\ell_{\alpha-1}$ through the isometries $(D_\alpha^{(j,+)}, d_\alpha) \equiv (D_{\alpha-1}, \frac{d_{\alpha-1}}{2})$ and $(D_\alpha^{(i,-)}, d_\alpha) \equiv (D_{\alpha-1}, \frac{d_{\alpha-1}}{2})$ (where $\alpha - 1$ denotes the predecessor of $\alpha$). To give meaning to the average of two elements, each belonging to one of these spaces, as an element of $\mathcal{F}(D_\alpha)$, we introduce the operators $L_\alpha^{(j,+)} \colon \mathcal{F}(D_\alpha^{(j,+)}) \to \mathcal{F}(D_\alpha)$ and $L_\alpha^{(i,-)} \colon \mathcal{F}(D_\alpha^{(i,-)}) \to \mathcal{F}(D_\alpha)$ defined by

$$L_\alpha^{(j,+)}(\delta_{D_\alpha^{(j,+)}}(x)) = \delta_{D_\alpha}(x) - \delta_{D_\alpha}(\ell_\alpha^{(j,+)})$$

and

$$L_\alpha^{(i,-)}(\delta_{D_\alpha^{(i,-)}}(x)) = \delta_{D_\alpha}(x) - \delta_{D_\alpha}(\ell_\alpha^{(i,-)})$$

such that they are linear.

**Fact 3.3.** *The maps $L_\alpha^{(j,+)}$ and $L_\alpha^{(i,-)}$ are well-defined isometries, satisfying*

$$\forall \gamma^+ \in \mathcal{F}(D_\alpha^{(j,+)}),\, \forall f \in \mathrm{Lip}_0(D_\alpha),\, \langle f,\, L_\alpha^{(j,+)}\gamma^+ \rangle = \langle f{\restriction}_{D_\alpha^{(j,+)}} - f(\ell_\alpha^{(j,+)}),\, \gamma^+ \rangle$$

*and*

$$\forall \gamma^- \in \mathcal{F}(D_\alpha^{(i,-)}),\, \forall f \in \mathrm{Lip}_0(D_\alpha),\, \langle f,\, L_\alpha^{(i,-)}\gamma^- \rangle = \langle f{\restriction}_{D_\alpha^{(i,-)}} - f(\ell_\alpha^{(i,-)}),\, \gamma^- \rangle.$$

*Proof.* The maps $L_\alpha^{(j,+)}$ and $L_\alpha^{(i,-)}$ can be seen as preadjoint operators. Indeed, let $Q\colon \mathrm{Lip}_0(D_\alpha) \to \mathrm{Lip}_0(D_\alpha^{(j,+)})$ defined by $Qf = f{\restriction}_{D_\alpha^{(j,+)}} - f(\ell_\alpha^{(j,+)})$. For every $f \in \mathrm{Lip}_0(D_\alpha)$, $\|Q(f)\|_L \leq \|f\|_L$ so $Q(B_{\mathrm{Lip}_0(D_\alpha)}) \subset B_{\mathrm{Lip}_0(D_\alpha^{(j,+)})}$. Conversely, given $f \in \mathrm{Lip}_0(D_\alpha^{(j,+)})$, it is easy to check that $g$ defined by $f$ on $D_\alpha^{(j,+)}$ and 0 at $\ell_\alpha$ is a Lipschitz map with $\|g\|_L \leq \|f\|_L$, so we can extend it to $g$ defined on $D_\alpha$ with $\|g\|_L \leq \|f\|_L$, and $Qg = f$. This proves that $B_{\mathrm{Lip}_0(D_\alpha^{(j,+)})} \subset Q(B_{\mathrm{Lip}_0(D_\alpha)})$. Since $Q$ is $w^*$-$w^*$-continuous, this implies the existence of an isometry $R\colon \mathcal{F}(D_\alpha^{(j,+)}) \to \mathcal{F}(D_\alpha)$ such that $Q = R^*$. Finally, for every $x \in D_\alpha^{(j,+)}$, for every $f \in \mathrm{Lip}_0(D_\alpha)$:

$$\begin{aligned}
\langle f,\, R\delta_{D_\alpha^{(j,+)}}(x) \rangle = \langle R^* f,\, \delta_{D_\alpha^{(j,+)}}(x) \rangle &= \langle Qf,\, \delta_{D_\alpha^{(j,+)}}(x) \rangle \\
&= \langle f{\restriction}_{D_\alpha^{(j,+)}} - f(\ell_\alpha^{(j,+)}),\, \delta_{D_\alpha^{(j,+)}}(x) \rangle \\
&= f(x) - f(\ell_\alpha^{(j,+)}) \\
&= \langle f,\, L_\alpha^{(j,+)} \delta_{D_\alpha^{(j,+)}}(x) \rangle,
\end{aligned}$$

so $R = L_\alpha^{(j,+)}$. We can do a similar discussion for $L_\alpha^{(i,-)}$. □

Now we have the following lemma:

**Lemma 3.4.** *Let $i \neq j \in \mathbb{N} \setminus \{1\}$, let $\alpha \in (0, \omega_1)$ be a successor ordinal, let $\varepsilon > 0$. Then, for every $\gamma^+ \in \mathcal{F}(D_\alpha^{(j,+)})$ and $\gamma^- \in \mathcal{F}(D_\alpha^{(i,-)})$:*

$$\min\left(\|\gamma^+\|_{\mathcal{F}(D_\alpha^{(j,+)})}, \|\gamma^-\|_{\mathcal{F}(D_\alpha^{(i,-)})}\right) \geq \varepsilon \Longrightarrow \left\| \frac{L_\alpha^{(j,+)}\gamma^+ + L_\alpha^{(i,-)}\gamma^-}{2} \right\|_{\mathcal{F}(D_\alpha)} \geq \varepsilon.$$

*Proof.* By assumption, there exists $f^+ \in \mathrm{Lip}_0(D_\alpha^{(j,+)})$ such that $\|f^+\|_L = 1$ and $\langle f^+,\ \gamma^+\rangle \geq \varepsilon$, and $f^- \in \mathrm{Lip}_0(D_\alpha^{(i,-)})$ such that $\|f^-\|_L = 1$ and $\langle f^-,\ \gamma^-\rangle \geq \varepsilon$. Let us define a function $f$ on $D_\alpha^{(j,+)} \cup D_\alpha^{(i,-)} \cup \{\ell_\alpha\}$ by $f = f^+$ on $D_\alpha^{(j,+)}$, $f = f^-$ on $D_\alpha^{(i,-)}$ and $f(\ell_\alpha) = 0$. Next, let us check that $f$ is 1-Lipschitz on $D_\alpha^{(j,+)} \cup D_\alpha^{(i,-)} \cup \{\ell_\alpha\}$. Notice that $f^+(\ell_\alpha^{(j,+)}) = f^-(\ell_\alpha^{(i,-)}) = 0$.

Consider first $x \in D_\alpha^{(j,+)}$ and $y \in D_\alpha^{(i,-)}$. Since $f^+$ and $f^-$ are 1-Lipschitz, we have

$$\begin{aligned}
|f(x) - f(y)| &= |f^+(x) - f^-(y)| \\
&\leq |f^+(x) - f^+(\ell_\alpha^{(j,+)})| + |f^-(\ell_\alpha^{(i,-)}) - f^-(y)| \\
&\leq \underbrace{d_\alpha(x, \ell_\alpha^{(j,+)})}_{\leq \frac{1}{2} + d_\alpha(x, t_\alpha)} + \underbrace{d_\alpha(\ell_\alpha^{(i,-)}, y)}_{\leq \frac{1}{2} + d_\alpha(x_\alpha^i, y)} \\
&\leq 1 + d_\alpha(x, t_\alpha) + d_\alpha(x_\alpha^i, y).
\end{aligned}$$

Note that $d_\alpha(x, \ell_\alpha^{(j,+)}) \leq \frac{1}{2} + d_\alpha(x, x_\alpha^j)$ and $d_\alpha(\ell_\alpha^{(i,-)}, y) \leq \frac{1}{2} + d_\alpha(b_\alpha, y)$, so we also have that $|f(x) - f(y)| \leq 1 + d_\alpha(x, x_\alpha^j) + d_\alpha(b_\alpha, y)$. Therefore,

$$|f(x) - f(y)| \leq 1 + \min\left(d_\alpha(x, t_\alpha) + d_\alpha(x_\alpha^i, y),\ d_\alpha(x, x_\alpha^j) + d_\alpha(b_\alpha, y)\right),$$

where the right-hand side is clearly equal to $d_\alpha(x, y)$.

Now for $x \in D_\alpha^{(j,+)}$,

$$\begin{aligned}
|f(x) - f(\ell_\alpha)| = |f^+(x)| &= |f^+(x) - f^+(\ell_\alpha^{(j,+)})| \\
&\leq d_\alpha(x, \ell_\alpha^{(j,+)}) \\
&\leq 1 + d_\alpha(x, t_\alpha)
\end{aligned}$$

where $1 + d_\alpha(x, t_\alpha)$ is equal to $d_\alpha(x, \ell_\alpha)$ since $x \in D_\alpha^{(j,+)}$ with $j \neq 1$. Similarly, for $x \in D_\alpha^{(i,-)}$, $|f(x) - f(\ell_\alpha)| \leq d_\alpha(x, \ell_\alpha)$. Thus $f$ is 1-Lipschitz on $D_\alpha^{(j,+)} \cup D_\alpha^{(i,-)} \cup \{\ell_\alpha\}$ so we can extend it to a 1-Lipschitz map defined on the whole $D_\alpha$, still denoted $f$. Thus, $f \in \mathrm{Lip}_0(D_\alpha)$. But according to Fact 3.3,

$$\begin{aligned}
\langle f,\ \frac{L_\alpha^{(j,+)}\gamma^+ + L_\alpha^{(i,-)}\gamma^-}{2}\rangle &= \frac{1}{2}\left(\langle f,\ L_\alpha^{(j,+)}\gamma^+\rangle + \langle f,\ L_\alpha^{(i,-)}\gamma^-\rangle\right) \\
&= \frac{1}{2}\left(\langle f^+,\ \gamma^+\rangle + \langle f^-,\ \gamma^-\rangle\right) \\
&\geq \frac{1}{2}(\varepsilon + \varepsilon) = \varepsilon.
\end{aligned}$$

since $f(\ell_\alpha^{(j,+)}) = f(\ell_\alpha^{(i,-)}) = 0$. So $\left\| \frac{L_\alpha^{(j,+)}\gamma^+ + L_\alpha^{(i,-)}\gamma^-}{2} \right\|_{\mathcal{F}(D_\alpha)} \geq \varepsilon$. □

**Proposition 3.5.** *Let $\alpha$, $\beta \in [0, \omega_1)$ such that $\alpha$ is a successor ordinal and $\alpha \geq \beta$. Let $i \neq j \in \mathbb{N} \setminus \{1\}$, let $\gamma^+ \in \mathcal{F}(D_\alpha^{(j,+)})$ and $\gamma^- \in \mathcal{F}(D_\alpha^{(i,-)})$. Assume that* $\gamma^+ \in \sigma_1^\beta(B_{\mathcal{F}(D_\alpha^{(j,+)})})$ *and* $\gamma^- \in \sigma_1^\beta(B_{\mathcal{F}(D_\alpha^{(i,-)})})$. *Then:*

$$\frac{L_\alpha^{(j,+)}\gamma^+ + L_\alpha^{(i,-)}\gamma^-}{2} \in \sigma_1^\beta(B_{\mathcal{F}(D_\alpha)}).$$

*Proof.* Let $\alpha$ be a successor ordinal. We proceed by transfinite induction on $\beta \leq \alpha$. The statement is immediate for $\beta = 0$. Assume now that it is true for every $\mu < \beta$. If $\beta$ is a limit ordinal, we get the conclusion by taking the intersection on $\mu < \beta$. If $\beta = \lambda + 1$ is a successor ordinal, let $\gamma^+ \in \sigma_1^{\lambda+1}(B_{\mathcal{F}(D_\alpha^{(j,+)})})$ and $\gamma^- \in \sigma_1^{\lambda+1}(B_{\mathcal{F}(D_\alpha^{(i,-)})})$. By induction hypothesis, $\frac{L_\alpha^{(j,+)}\gamma^+ + L_\alpha^{(i,-)}\gamma^-}{2} \in \sigma_1^\lambda(B_{\mathcal{F}(D_\alpha)})$. Let $V$ be a relative $w$-open neighborhood of $\frac{L_\alpha^{(j,+)}\gamma^+ + L_\alpha^{(i,-)}\gamma^-}{2}$ in $\sigma_1^\lambda(B_{\mathcal{F}(D_\alpha)})$. We must show that $\operatorname{diam}(V) \geq 1$. Without loss of generality, we may assume that $V$ is of the form

$$V = \{\gamma \in \sigma_1^\lambda(B_{\mathcal{F}(D_\alpha)}) : \forall r \in \{1, \dots, n\},\ \Big|\langle f_r,\ \gamma - \frac{L_\alpha^{(j,+)}\gamma^+ + L_\alpha^{(i,-)}\gamma^-}{2}\rangle\Big| \leq \varepsilon\}$$

with $\varepsilon > 0$ and $f_1, \dots, f_n \in \mathrm{Lip}_0(D_\alpha)$. By induction hypothesis and triangle inequality, $V$ contains the set $\frac{1}{2}(L_\alpha^{(j,+)}W^+ + L_\alpha^{(i,-)}W^-)$ where

$$\begin{aligned} W^+ &:= \{\gamma \in \sigma_1^\lambda(B_{\mathcal{F}(D_\alpha^{(j,+)})}) : \forall r \in \{1, \dots, n\},\ \big|\langle f_r,\ L_\alpha^{(j,+)}\gamma - L_\alpha^{(j,+)}\gamma^+\rangle\big| \leq \varepsilon\} \\ &\overset{\text{Fact 3.3}}{=} \Big\{\gamma \in \sigma_1^\lambda(B_{\mathcal{F}(D_\alpha^{(j,+)})}) : \forall r \in \{1, \dots, n\}, \\ &\qquad \big|\langle f_r\restriction_{D_\alpha^{(j,+)}} - f_r(\ell_\alpha^{(j,+)}),\ \gamma - \gamma^+\rangle\big| \leq \varepsilon\Big\} \end{aligned}$$

is a relative $w$-open neighborhood of $\gamma^+$ in $\sigma_1^\lambda(B_{\mathcal{F}(D_\alpha^{(j,+)})})$, and

$$\begin{aligned} W^- &:= \{\gamma \in \sigma_1^\lambda(B_{\mathcal{F}(D_\alpha^{(i,-)})}) : \forall r \in \{1, \dots, n\},\ \big|\langle f_r,\ L_\alpha^{(i,-)}\gamma - L_\alpha^{(i,-)}\gamma^-\rangle\big| \leq \varepsilon\} \\ &\overset{\text{Fact 3.3}}{=} \Big\{\gamma \in \sigma_1^\lambda(B_{\mathcal{F}(D_\alpha^{(i,-)})}) : \forall r \in \{1, \dots, n\}, \\ &\qquad \big|\langle f_r\restriction_{D_\alpha^{(i,-)}} - f_r(\ell_\alpha^{(i,-)}),\ \gamma - \gamma^-\rangle\big| \leq \varepsilon\Big\} \end{aligned}$$

is a relative $w$-open neighborhood of $\gamma^-$ in $\sigma_1^\lambda(B_{\mathcal{F}(D_\alpha^{(i,-)})})$. Thus, $\operatorname{diam}(W^+) \geq 1$ and $\operatorname{diam}(W^-) \geq 1$ so, for every $\eta < 1$, there exist $\mu^+, \nu^+ \in W^+$ and $\mu^-, \nu^- \in W^-$ such that $\|\mu^+ - \nu^+\| \geq \eta$ and $\|\mu^- - \nu^-\| \geq \eta$. Setting $u := \frac{L_\alpha^{(j,+)}\mu^+ + L_\alpha^{(i,-)}\mu^-}{2}$ and $v := \frac{L_\alpha^{(j,+)}\nu^+ + L_\alpha^{(i,-)}\nu^-}{2}$, $u$ and $v$ are two elements of $\frac{1}{2}(L_\alpha^{(j,+)}W^+ + L_\alpha^{(i,-)}W^-) \subset V$ satisfying

$$\|u - v\| = \left\|\frac{L_\alpha^{(j,+)}(\mu^+ - \nu^+) + L_\alpha^{(i,-)}(\mu^- - \nu^-)}{2}\right\| \geq \eta$$

according to Lemma 3.4. Therefore, $\operatorname{diam}(V) \geq \eta$ for all $\eta < 1$, which concludes the proof. □

3.3. **Higher derived sets.** We are now able to prove the main result of this section:

**Proposition 3.6.** *For every* $\alpha \in (0, \omega_1)$, $m_{t_\alpha, b_\alpha} \in \sigma_1^\alpha(B_{\mathcal{F}(D_\alpha)})$.

*Proof.* We will use a transfinite induction on $\alpha$. The statement is true for $\alpha = 1$ thanks to Proposition 3.1. Assume it is true for every $\beta < \alpha$. If $\alpha$ is a limit ordinal, for every $\beta < \alpha$, $m_{t_\beta,b_\beta} \in \sigma_1^\beta(B_{\mathcal{F}(D_\beta)})$ by induction hypothesis. Setting $U_\beta \colon \mathcal{F}(D_\beta) \to \mathcal{F}(D_\alpha)$ the isometry such that $U_\beta(\delta_{D_\beta}(x)) = \delta_{D_\alpha}((\beta,x)) - \delta_{D_\alpha}((\beta,\ell_\beta))$, this implies that $U_\beta(m_{t_\beta,b_\beta}) = m_{t_\alpha,b_\alpha}$ belongs to $\sigma_1^\beta(B_{\mathcal{F}(D_\alpha)})$. Since this holds for every $\beta < \alpha$, we have $m_{t_\alpha,b_\alpha} \in \sigma_1^\alpha(B_{\mathcal{F}(D_\alpha)})$. So assume now $\alpha = \beta+1$, and let $V$ be a $w$-open neighborhood of $m_{t_\alpha,b_\alpha}$ in $\mathcal{F}(D_\alpha)$. Recall that according to Remark 3.2, there exist $j > i$ in $\mathbb{N} \setminus \{1\}$ such that $\gamma_V \coloneqq \frac{1}{2}(m_{t_\alpha,x_\alpha^j} + m_{x_\alpha^i,b_\alpha}) \in V$. However, by definition of $D_\alpha$, $m_{t_\alpha,x_\alpha^j}$ identifies with $m_{t_\beta,b_\beta} \in \mathcal{F}(D_\beta,d_\beta) \cong \mathcal{F}(D_\alpha^{(j,+)},2d_\alpha)$, so $\frac{\delta(t_\alpha)-\delta(x_\alpha^j)}{2} \in \sigma_1^\beta(B_{\mathcal{F}(D_\alpha^{(j,+)},\,2d_\alpha)})$ by induction hypothesis. Setting $U \colon \mathcal{F}(D_\alpha^{(j,+)},\,2d_\alpha) \to \mathcal{F}(D_\alpha^{(j,+)},d_\alpha)$ the linear operator such that $U(\delta(x)) = 2\delta(x)$, this implies that $U\left(\frac{\delta(t_\alpha)-\delta(x_\alpha^j)}{2}\right) \in \sigma_1^\beta(B_{\mathcal{F}(D_\alpha^{(j,+)})})$ because $U$ is an isometry, that is, $m_{t_\alpha,x_\alpha^j} \in \sigma_1^\beta(B_{\mathcal{F}(D_\alpha^{(j,+)})})$. Similarly, $m_{x_\alpha^i,b_\alpha} \in \sigma_1^\beta(B_{\mathcal{F}(D_\alpha^{(i,-)})})$, so Proposition 3.5 yields $\frac{L_\alpha^{(j,+)}(m_{t_\alpha,x_\alpha^j})+L_\alpha^{(i,-)}(m_{x_\alpha^i,b_\alpha})}{2} = \gamma_V \in \sigma_1^\beta(B_{\mathcal{F}(D_\alpha)})$. Observe that the net $(\gamma_V)_V \subset \sigma_1^\beta(B_{\mathcal{F}(D_\alpha)})$ is $w$-convergent to $m_{t_\alpha,b_\alpha}$, with $\sigma_1^\beta(B_{\mathcal{F}(D_\alpha)})$ $w$-closed: hence, $m_{t_\alpha,b_\alpha} \in \sigma_1^\beta(B_{\mathcal{F}(D_\alpha)})$. Since $\gamma_V$ belongs to $V \cap \sigma_1^\beta(B_{\mathcal{F}(D_\alpha)})$ and $\|\gamma_V - m_{t_\alpha,b_\alpha}\| = 1$, it follows that $\operatorname{diam}(V \cap \sigma_1^\beta(B_{\mathcal{F}(D_\alpha)})) \geq 1$ and thus $m_{t_\alpha,b_\alpha} \in \sigma_1^{\beta+1}(B_{\mathcal{F}(D_\alpha)})$. □

In particular, $\sigma_1^\alpha(B_{\mathcal{F}(D_\alpha)}) \neq \varnothing$ so $\Phi(\mathcal{F}(D_\alpha),1) > \alpha$. Finally, we obtain that there exist Lipschitz-free spaces satisfying the PCP "as badly as possible" (*cf* Proposition 2.5):

**Theorem 3.7.** *For every $\alpha \in (0,\omega_1)$, there exists a countable complete metric space $D_\alpha$ such that $\Phi(\mathcal{F}(D_\alpha)) > \alpha$.*

**Remark 3.8.** We can draw a parallel between this theorem and a result of Braga, Lancien, Petitjean and Procházka. Indeed, in [6] (Theorem 4.3), they exhibited a uniformly discrete metric space $M$ such that for each Banach space whose dual contains an isomorphic copy of $\mathcal{F}(M)$, the Szlenk index of this space is greater than $\omega^2$. Here, notice that for every Banach space $X$, the Szlenk index of $X$ is greater than $\Phi(X^*)$. Then we have that for each Banach space $X$ whose dual contains an isomorphic copy of $\mathcal{F}(D_\alpha)$, the Szlenk index of $X$ is greater than $\alpha$.

3.4. **Computation for $\alpha = \omega$.** Now that we have a lower bound for $\Phi(\mathcal{F}(D_\alpha))$, a natural question is whether we can compute its exact value. In order to do that, we can try to apply the following result, which is well-known to specialists. Until the end of this section, $(X_n)_{n\in\mathbb{N}}$ will stand for a family of Banach spaces, and we set $X = (\sum_{n\in\mathbb{N}} X_n)_{\ell_1}$. Then:

**Proposition 3.9.** *For every $\varepsilon \in (0,1)$, we have:*

$$\Phi(X,3\varepsilon) \leq \sup_{n\in\mathbb{N}} \alpha_n \,.\, \omega,$$

*writing $\alpha_n \coloneqq \max_{1\leq k\leq n} \Phi(X_k)$.*

To show this proposition, we can use for example the following lemma, which is an adaptation of Lemma 3.3 in [13]:

**Lemma 3.10.** *Let $\varepsilon \in (0,1)$, let $z \in B_X$ and let $n \in \mathbb{N}$ such that $\|P_n z\| > 1-\varepsilon$, where $P_n$ denotes the canonical projection from $X$ onto $(\sum_{k=1}^{n} X_k)_{\ell_1}$. Then, for every ordinal $\alpha \in [0,\, \omega_1)$:*

$$z \in \sigma_{3\varepsilon}^{\alpha}(B_X) \implies P_n z \in \sigma_{\varepsilon}^{\alpha}(P_n B_X).$$

*Proof.* We will use a transfinite induction on $\alpha$. The statement is clearly true for $\alpha = 0$. If it is true for every $\beta < \alpha$, then if $\alpha$ is a limit ordinal, it is also true for $\alpha$ by taking the intersection on $\beta < \alpha$. So assume now that $\alpha = \mu + 1$, and let $z \in B_X$ and $n \in \mathbb{N}$ such that $\|P_n z\| > 1-\varepsilon$. We proceed by contraposition: assume that $P_n z \notin \sigma_{\varepsilon}^{\alpha}(P_n B_X)$, and let us show that $z \notin \sigma_{3\varepsilon}^{\alpha}(B_X) = \sigma_{3\varepsilon}(\sigma_{3\varepsilon}^{\mu}(B_X))$. So we may also assume that $z \in \sigma_{3\varepsilon}^{\mu}(B_X)$. Then the induction hypothesis implies that $P_n z \in \sigma_{\varepsilon}^{\mu}(P_n B_X)$, and since $P_n z \notin \sigma_{\varepsilon}^{\alpha}(P_n B_X) = \sigma_{\varepsilon}(\sigma_{\varepsilon}^{\mu}(P_n B_X))$, there exists $V$ a $w$-open subset of $P_n X = \sum_{k=1}^{n} X_k$ containing $P_n z$ such that $\operatorname{diam}(V \cap \sigma_{\varepsilon}^{\mu}(P_n B_X)) < \varepsilon$. We may assume that $V$ is of the form

$$V = \Big\{x \in \sum_{k=1}^{n} X_k : \forall i \in \{1,\dots,r\},\ f_i(x) > \alpha_i\Big\}$$

with $\alpha_i \in \mathbb{R}$ and $f_i \in (\sum_{k=1}^{n} X_k)_{\ell_1}^{*}$ of norm one. Since $\|P_n z\| > 1-\varepsilon$, we may also assume that $\alpha_1 > 1-\varepsilon$. This last assumption implies that $V \cap (1-\varepsilon)B_X = \varnothing$. Now we extend each $f_i$ to $g_i \in X^*$ by setting $g_i = f_i$ on $\sum_{k=1}^{n} X_k$ and $g_i = 0$ on $\sum_{k>n} X_k$. Setting

$$U = \{x \in X : \forall i \in \{1,\dots,r\},\ g_i(x) > \alpha_i\},$$

we can notice that $U$ is a $w$-open subset of $X$ containing $z$, so $z \in U \cap \sigma_{3\varepsilon}^{\mu}(B_X)$. To conclude that $z \notin \sigma_{3\varepsilon}^{\mu+1}(B_X)$, it remains to show that $\operatorname{diam}(U \cap \sigma_{3\varepsilon}^{\mu}(B_X)) < 3\varepsilon$: let $x, y \in U \cap \sigma_{3\varepsilon}^{\mu}(B_X)$. From the definition of the $g_i$'s, we have that $P_n x$, $P_n y$ belong to $V$ and thus are of norm strictly larger than $1-\varepsilon$. Since $\|x\| = \|P_n x\| + \|x - P_n x\|$, it follows that $\|x - P_n x\| < \varepsilon$, and likewise $\|y - P_n y\| < \varepsilon$. So

$$\begin{aligned}\|x-y\| &\le \|x - P_n x\| + \|P_n x - P_n y\| + \|P_n y - y\| \\ &\le 2\varepsilon + \operatorname{diam}(V \cap \sigma_{\varepsilon}^{\mu}(P_n B_X)) < 3\varepsilon.\end{aligned}$$

Therefore, $\operatorname{diam}(U \cap \sigma_{3\varepsilon}^{\mu}(B_X)) < 3\varepsilon$. □

Now we have all the tools to prove the desired proposition:

*Proof of Proposition 3.9.* Let $x \in B_X$ such that $\|x\| > 1-\varepsilon$ and let $n \in \mathbb{N}$ such that $\|P_n x\| > 1-\varepsilon$. We will use the following known fact:

**Fact 3.11.** *For every $n \in \mathbb{N}$, we have $\Phi((\sum_{k=1}^{n} X_k)_{\ell_1}, \varepsilon) \le \max_{1\le k\le n} \Phi(X_k)$.*

*Ideas of the proof of Fact 3.11.* The first point is that for every Banach spaces $X$, $Y$, for every $\varepsilon > 0$, $\Phi(X \oplus_\infty Y, \varepsilon) \le \max(\Phi(X), \Phi(Y))$. Indeed, it is proved in Proposition 2.4 in [13] that the Szlenk index of $X \oplus_1 X$ is equal to the Szlenk index of $X$. The authors did the proof considering diameters measured in the sum $X^* \oplus_\infty X^*$, here it is a slight modification considering diameters measured in the sum $X \oplus_\infty Y$. We would like to draw the reader's attention to a minor oversight in an inclusion used in the proof of Proposition 2.4 in [13]: the inclusion $\sigma_{\varepsilon}^{\alpha}(C \cup D) \subset \sigma_{\varepsilon}^{\alpha}(C) \cup \sigma_{\varepsilon}^{\alpha}(D)$ is not true in general. Instead, the correct inclusion is $\sigma^{\alpha}_{\varepsilon}(C \cup D) \subset \sigma^{\alpha}_{\varepsilon/2}(C) \cup \sigma^{\alpha}_{\varepsilon/2}(D)$; see for instance Lemma 3.1 in [7] for a proof of this fact. However, this oversight does not affect the desired result in any way.

Since the property "for every $\varepsilon > 0$, $\Phi(X \oplus_\infty Y, \varepsilon) \leq \max(\Phi(X), \Phi(Y))$" is an isomorphic property, we now have the following inequality:

$$\forall \varepsilon > 0,\ \Phi(X \oplus_1 Y, \varepsilon) \leq \max(\Phi(X), \Phi(Y)).$$

Thus, $\Phi(X \oplus_1 Y) \leq \max(\Phi(X), \Phi(Y))$, and an iterated application of this leads to Fact 3.11. ☐

Let us get back to the proof of Proposition 3.9. Using the previous fact, we have $\sigma^{\alpha_n}_{\varepsilon}(P_n B_X) = \varnothing$ and then Lemma 3.10 implies that $x \notin \sigma^{\alpha_n}_{3\varepsilon}(B_X)$. Setting $\alpha := \sup_{n \in \mathbb{N}} \alpha_n$, the previous sentence implies $\sigma^{\alpha}_{3\varepsilon}(B_X) \subset (1-\varepsilon)B_X$. A homogeneity argument leads to

$$\forall k \in \mathbb{N},\ \sigma^{\alpha.k}_{3\varepsilon}(B_X) \subset (1-\varepsilon)^k B_X.$$

Considering $k \in \mathbb{N}$ such that $(1-\varepsilon)^k < \frac{3\varepsilon}{2}$, we have that $(1-\varepsilon)^k B_X$ is of diameter strictly less than $3\varepsilon$, and thus $\sigma^{\alpha.(k+1)}_{3\varepsilon}(B_X) = \varnothing$, which concludes the proof. ☐

**Remark 3.12.** In particular, combining Kalton's decomposition (see Proposition 4.3 in [16]) with Proposition 3.9, we obtain what seems to be a folklore fact among specialists: for every uniformly discrete metric space $M$, $\Phi(\mathcal{F}(M)) \leq \omega^2$ (indeed, the free space over a bounded uniformly discrete space is isomorphic to $\ell_1$, and $\Phi(\ell_1) = \omega$).

Here, $D_\omega$ is not uniformly discrete and not even discrete. However, we are still able to show the following:

**Proposition 3.13.** *The index $\Phi(\mathcal{F}(D_\omega))$ is equal to $\omega^2$.*

*Proof.* Let us consider the open covering of $D_\omega$ given by the sets

$$A := \{z \in D_\omega : d_\omega(z, b_\omega) < \frac{3}{2}\} \text{ and } B := \{z \in D_\omega : d_\omega(z, t_\omega) < \frac{3}{2}\},$$

and the function defined by $D(z) := \mathrm{dist}(z, D_\omega \setminus A) + \mathrm{dist}(z, D_\omega \setminus B)$ for $z \in D_\omega$. We wish to apply Lemma 2.5 in [1] with this covering. To this end, we must check that $\inf_{z \in D_\omega} D(z) > 0$. There are three cases:

*a*) If $z \in A \setminus B$, we have $D(z) = \mathrm{dist}(z, D_\omega \setminus A)$. Either $z = b_\omega$ and then $D(z) \geq \frac{3}{2}$, or $z = (n, x)$ for some $n \in \mathbb{N}$ and $x \in D_n$, so the closest point to $z$ in $D_\omega \setminus A$ is of the form $(n, y)$ with $y \in D_n$ and $d_n(y, b_\omega) \geq \frac{3}{2}$. Then

$$D(z) = d_n(x, y) = d_n(y, b_\omega) - d_n(x, b_\omega) \geq \frac{3}{2} - \frac{1}{2} = 1.$$

*b*) If $z \in B \setminus A$, by symmetry with the first case we have $D(z) \geq 1$.

*c*) If $z \in A \cap B$, we have $D(z) = \mathrm{dist}(z, D_\omega \setminus A) + \mathrm{dist}(z, D_\omega \setminus B)$. Either $d_\omega(z, b_\omega) \leq 1$ and then $\mathrm{dist}(z, D_\omega \setminus A) \geq \frac{1}{2}$, or $d_\omega(z, t_\omega) \leq 1$ and then $\mathrm{dist}(z, D_\omega \setminus B) \geq \frac{1}{2}$. Thus, $D(z) \geq \frac{1}{2}$.

Consequently, $\inf_{z \in D_\omega} D(z) \geq \frac{1}{2} > 0$, so we can apply Lemma 2.5 followed by Lemma 2.4 in [1] to obtain that $\mathcal{F}(D_\omega)$ is isomorphic to a subspace of $\mathcal{F}(A) \oplus \mathcal{F}(B)$.

Now, for $n \in \mathbb{N}$, let

$$A_n := \{z \in D_\omega : z = (n, x) \text{ with } x \in D_n \text{ and } d_\omega(z, b_\omega) < \frac{3}{2}\}.$$

In order to get an upper estimate for $\Phi(\mathcal{F}(D_\omega))$, it is enough to get upper estimates for $\Phi(\mathcal{F}(A))$ and $\Phi(\mathcal{F}(B))$. For this purpose, we will use a result in [26] to write $\mathcal{F}(A)$ as the $\ell_1$-sum of the $\mathcal{F}(A_n)$. Without loss of generality, take the base point of $A_n$ and $A$ to be $b_\omega$, so that $A = \dot{\bigcup\limits_{n\in\mathbb{N}}} A_n$. Let us denote by $d^1$ the summing distance on $A$: $d^1$ is defined by $d^1\restriction_{A_n\times A_n} = d_n$ and $d^1(x,y) = d_\omega(x,b_\omega)+d_\omega(b_\omega,y)$ whenever $x$ and $y$ belong to distinct summands. According to Proposition 3.9 in [26],

$$\mathcal{F}(A,d^1) \cong (\sum_{n\in\mathbb{N}} \mathcal{F}(A_n))_{\ell_1}.$$

But $d^1$ and $d_\omega$ are Lipschitz equivalent on $A$; more precisely, $\frac{1}{3}d^1 \leq d_\omega \leq d^1$. Indeed, if $x$ and $y$ are in distinct summands and if $d_\omega(x,y) \neq d^1(x,y)$, then $d_\omega(x,y) = d_\omega(x,t_\omega) + d_\omega(y,t_\omega) > \frac{1}{2} + \frac{1}{2} = 1$ by definition of $A$, and thus $d^1(x,y) < \frac{3}{2} + \frac{3}{2} = 3 < 3d_\omega(x,y)$. The other inequality results directly from the definition of $d_\omega$. As a consequence,

$$\mathcal{F}(A,d_\omega) \simeq \mathcal{F}(A,d^1) \cong (\sum_{n\in\mathbb{N}} \mathcal{F}(A_n))_{\ell_1}.$$

Using Proposition 3.9, we can deduce that $\Phi(\mathcal{F}(A,d_\omega)) \leq \omega^2$, since the $A_n$ are bounded and uniformly discrete.

Applying the same reasoning to $B$ yields

$$\Phi(\mathcal{F}(D_\omega)) \leq \max\left(\Phi(\mathcal{F}(A,d_\omega)), \Phi(\mathcal{F}(B,d_\omega))\right) \leq \omega^2.$$

Finally, as $\Phi(\mathcal{F}(D_\omega)) > \omega$, Proposition 2.2 leads to $\Phi(\mathcal{F}(D_\omega)) = \omega^2$. □

**Remark 3.14.** Drawing inspiration from the computation of the Szlenk index of $\mathcal{C}(K)$ spaces in [24], it is tempting to conjecture that given $\beta \in [0,\omega_1)$, for all $\alpha \in [\omega^\beta,\ \omega^{\beta+1})$, we have $\Phi(\mathcal{F}(D_\alpha)) = \omega^{\beta+1}$.

It is possible to adjust the proof above to show that $\Phi(\mathcal{F}(D_\alpha)) \leq \sup\limits_{\beta<\alpha} \Phi(\mathcal{F}(D_\beta)) \,.\, \omega$ for every limit ordinal $\alpha < \omega_1$. Adapting Theorem 2.11 in [7], it is also possible to show that for every $\alpha < \omega_1$, $\Phi(\mathcal{F}(D_\alpha)) = \Phi(\mathcal{F}(D_{\alpha+1}))$. But this yields an upper estimate for $\Phi(\mathcal{F}(D_\alpha))$ too rough to prove our conjecture.

## 4. Consequences

### 4.1. **Universal spaces.**

If $\mathcal{C}$ is a class of metric spaces, we say that a metric space $M$ is *Lipschitz-universal* for the class $\mathcal{C}$ if every member of $\mathcal{C}$ Lipschitz-embeds into $M$, that is, is Lipschitz-equivalent to a subspace of $M$.

Using Theorem 3.7, we can deduce that a separable complete Lipschitz-universal space for the class of countable complete metric spaces cannot be purely 1-unrectifiable:

**Corollary 4.1.** *Let $M$ be a separable complete metric space such that every countable complete metric space is Lipschitz-equivalent to a subspace of $M$. Then $M$ is not purely 1-unrectifiable.*

*Proof.* For every $\alpha \in (0,\omega_1)$, by Proposition 2.7, $D_\alpha$ Lipschitz-embeds into $M$, so $\mathcal{F}(D_\alpha)$ is linearly isomorphic to a subspace of $\mathcal{F}(M)$. Consequently, for every $\alpha \in (0,\omega_1)$ we have $\Phi(\mathcal{F}(M)) > \alpha$, so $\Phi(\mathcal{F}(M)) \geq \omega_1$. Using Proposition 2.5, we deduce that $\mathcal{F}(M)$ does not have the PCP, which finishes the proof, given the characterization recalled in Remark 1.1. □

In particular, a separable complete metric space which is Lipschitz-universal for purely 1-unrectifiable complete metric spaces cannot be purely 1-unrectifiable: in a sense, there does not exist a separable Lipschitz-universal space for purely 1-unrectifiable complete metric spaces.

4.2. **Non-isomorphic Lipschitz-free spaces over purely 1-unrectifiable metric spaces.** For each $\alpha \in (0, \omega_1)$, according to Proposition 2.5 and Theorem 3.7, we have that $\alpha < \Phi(\mathcal{F}(D_\alpha)) < \omega_1$. Then, with a transfinite induction, it is easy to build a map $\varphi\colon (0, \omega_1) \to (0, \omega_1)$ such that:

$$\forall \alpha < \beta,\ \Phi(\mathcal{F}(D_{\varphi(\beta)})) > \Phi(\mathcal{F}(D_{\varphi(\alpha)})).$$

Since the weak-fragmentability index is an isomorphic invariant, we deduce the following result:

**Corollary 4.2.** *There exists an uncountable family $(M_i)_{i\in I}$ of countable complete metric spaces such that their Lipschitz-free spaces $(\mathcal{F}(M_i))_{i\in I}$ are pairwise non-isomorphic.*

That there are uncountably many non-isomorphic Lipschitz-free spaces over separable metric spaces was proved for the first time by Hájek, Lancien and Pernecká in [14] using a very different method. Indeed, their family consists of free spaces over separable Banach spaces.

4.3. **Lipschitz-free spaces over a compact.** It is still an open question whether for every separable Banach space $X$, there exists a compact space $K$ such that $\mathcal{F}(X)$ is linearly isomorphic to $\mathcal{F}(K)$. For example, the answer is positive for finite-dimensional spaces (Corollary 3.3 in [17] states that $\mathcal{F}(X)$ is linearly isomorphic to $\mathcal{F}(B_X)$) and for the Pełczyński universal space $\mathbb{P}$ (see [11]). This question was also open when considering a metric space $M$ instead of $X$ (see Question 1.2); the following corollary provides a negative answer to it:

**Corollary 4.3.** *Let $\alpha \in [\omega, \omega_1)$ and $K$ be any compact metric space. Then $\mathcal{F}(D_\alpha)$ and $\mathcal{F}(K)$ are not isomorphic.*

Let us recall some notions which will be involved in the proof: if $X$ is a Banach space and $S_X$ is its unit sphere, the *modulus of asymptotic uniform convexity* of X is given by $\bar{\delta}_X(t) = \inf_{x\in S_X} \bar{\delta}_X(t,x)$ for $t > 0$, where

$$\bar{\delta}_X(t,x) = \sup_{\dim(X/Y)<\infty} \inf_{y\in S_Y} \|x+ty\| - 1$$

and its *modulus of asymptotic uniform smoothness* is given by $\bar{\rho}_X(t) = \sup_{x\in S_X} \bar{\rho}_X(t,x)$ for $t > 0$, where

$$\bar{\rho}_X(t,x) = \inf_{\dim(X/Y)<\infty} \sup_{y\in S_Y} \|x+ty\| - 1.$$

We say that $X$ is *asymptotically uniformly convex* (AUC for short) if $\bar{\delta}_X(t) > 0$ for every $t > 0$ and that $X$ is *asymptotically uniformly smooth* (AUS for short) if $\lim_{t\to 0} t^{-1}\bar{\rho}_X(t) = 0$. If $X$ is a dual space and if we consider $w^*$-closed subspaces $Y$ of $X$ instead of norm-closed subspaces, then we denote $\bar{\delta}^*_X(t)$ the corresponding modulus and we say that $X$ is *weak* asymptotically uniformly convex* (AUC* for short) if $\bar{\delta}^*_X(t) > 0$ for every $t > 0$. It is well known (see Proposition 2.8 in [12] for the separable case and Corollary 2.4 in [10] for the general case) that a Banach space is AUS if and only if its dual space is AUC*. Finally, we say that $X$ is *AUC* *renormable* (resp. *AUS renormable*) if $X$ admits an equivalent AUC (resp. AUS) norm.

It is also well known that the notion of asymptotic uniform convexity is related to the weak-fragmentability index in the following way:

**Proposition 4.4.** *Let $X$ be a Banach space. If $X$ is AUC renormable, then $\Phi(X) \leq \omega$.*

*Proof.* Let $\varepsilon \in (0,1)$. It is easy to check that

$$\sigma_\varepsilon(B_X) \subset (1-\Delta(\varepsilon))B_X$$

where $\Delta(\varepsilon) = \frac{1}{2}\overline{\delta}_X(\frac{\varepsilon}{3})$. Then, an homogeneity argument leads to

$$\sigma_\varepsilon^n(B_X) \subset (1-\Delta(\varepsilon))^n B_X$$

for every $n \in \mathbb{N}$. Let $n \in \mathbb{N}$ be large enough so that $(1-\Delta(\varepsilon))^n < \frac{\varepsilon}{2}$. Then, since $\frac{\varepsilon}{2}B_X$ is of diameter strictly less than $\varepsilon$, we have that $\sigma_\varepsilon^{n+1}(B_X) = \varnothing$, which concludes the proof. □

Now we have all the tools to prove the corollary:

*Proof of Corollary 4.3.* Assume that there are $\alpha \in [\omega, \omega_1)$ and some compact metric space $K$ such that $\mathcal{F}(D_\alpha) \simeq \mathcal{F}(K)$. Since $D_\alpha$ is purely 1-unrectifiable, $\mathcal{F}(D_\alpha)$ has the PCP, and so does $\mathcal{F}(K)$. Using again Remark 1.1, $K$ is a purely 1-unrectifiable compact space. Thus, with Theorem 3.2 in [2], $\mathcal{F}(K)$ is isometric to the dual of $\mathrm{lip}_0(K)$, the space of locally flat Lipschitz functions on $K$ vanishing at 0. But Kalton proved in [16] (Theorem 6.6) that whenever $K$ is compact, $\mathrm{lip}_0(K)$ is isomorphic to a subspace of $c_0$. Since $c_0$ is AUS, it follows that $\mathrm{lip}_0(K)$ is AUS renormable. So $\mathcal{F}(K) \equiv \mathrm{lip}_0(K)^*$ is AUC$^*$ renormable, and hence AUC renormable. As a consequence of Proposition 4.4 we have $\Phi(\mathcal{F}(K)) \leq \omega$ and thus $\Phi(\mathcal{F}(D_\alpha)) \leq \omega$, a contradiction. □

**Remark 4.5.** Actually, in [16] (Theorem 6.6), Kalton proved a more precise result: if $K$ is a compact metric space, then for every $\varepsilon > 0$, $\mathrm{lip}_0(K)$ is $(1+\varepsilon)$-isomorphic to a subspace of $c_0$. Therefore, Lemma 4.4.1 in [22] shows that $\mathrm{lip}_0(K)$ is AUS, and not only AUS renormable (and thus $\mathcal{F}(K)$ is AUC, and not only AUC renormable). Showing that $\mathcal{F}(K)$ is AUC renormable is enough to get a contradiction in the previous proof, but this observation will be useful in the sequel (see Proposition 4.7).

4.4. **Compact reduction.** Proposition 4.4 implies that the space $\mathcal{F}(D_\omega)$ is not AUC renormable. But actually, this fact could also be deduced from the next theorem which was proved in [4] (Corollary 5.3):

**Theorem 4.6.** *If the family of the countably branching diamond graphs $(D_k)_{k\in\mathbb{N}}$ equi-Lipschitz embeds into a Banach space $X$, then $X$ is not AUC renormable.*

However, we can still give an interesting fact about the property "being AUC". A Banach space property $\mathcal{P}$ is said to be *compactly determined* if a Lipschitz-free space $\mathcal{F}(M)$ has $\mathcal{P}$ whenever the subspace $\mathcal{F}(K)$ has $\mathcal{P}$ for each compact $K \subset M$. For example, the weak sequential completeness or the Schur property are compactly determined properties. See [3] for more information on this subject.

With Corollary 4.3, we deduce:

**Proposition 4.7.** *The following Banach space properties are not compactly determined:*

*(i)* *being AUC;*
*(ii)* *being AUC renormable;*
*(iii)* *having a weak-fragmentability index lower than* $\omega$*;*
*(iv)* *having a weak-fragmentability index lower than* $\beta$ *for some fixed* $\beta \in (0, \omega_1)$*.*

*Proof.* If $K$ is a compact subset of $D_\alpha$ for some $\alpha \in (0, \omega_1)$, then $K$ is a countable compact metric space so a result of Dalet gives that $\mathcal{F}(K)$ is isometric to $\mathrm{lip}_0(K)^*$ (see Theorem 2.1 in [9]). As in the proof of Corollary 4.3, we deduce that $\mathcal{F}(K)$ is AUC renormable. Therefore, $\mathcal{F}(K)$ is AUC renormable for each compact $K \subset D_\alpha$ (and thus $\Phi(\mathcal{F}(K)) \leq \omega$), while $\mathcal{F}(D_\alpha)$ is not (and $\Phi(\mathcal{F}(D_\alpha)) > \omega$ for $\alpha \in (\omega, \omega_1)$), so *(ii)* and *(iii)* are not compactly determined properties.

With Remark 4.5, we adapt the above to obtain *(i)*.

Finally, given $\beta \in (0, \omega_1)$, we consider two cases. First, if $\beta \geq \omega$, we have $\Phi(\mathcal{F}(D_\beta)) > \beta$ while as above, $\Phi(\mathcal{F}(K)) \leq \omega \leq \beta$ for every compact subset $K$ of $D_\beta$, since subsets of $D_\beta$ are countable. On the other hand, if $\beta = n \in \mathbb{N}$, $D_\beta$ is uniformly discrete so a compact subset $K$ of $D_\beta$ must be a finite set. Then the weak topology and the norm topology coincide on $\mathcal{F}(K)$ because $\mathcal{F}(K)$ is finite-dimensional, so $\Phi(\mathcal{F}(K)) = 1 \leq \beta$, while $\Phi(\mathcal{F}(D_\beta)) > \beta$. □

## 5. Final remarks

5.1. **Corollaries 4.1 and 4.2 using dentability.** Considering slices instead of weakly open sets in the definition of the weak-fragmentability index, we obtain the dentability index $D(X)$ of a Banach space $X$; see [20] for a more detailed definition and some elementary facts. Like $\Phi$, the index $D$ has the following properties: it is an isomorphic invariant, if $F \subset X$ is a closed subspace of $X$ then $D(F) \leq D(X)$, and if $X$ is a separable Banach space with the RNP, then $D(X) < \omega_1$.

A result similar to Theorem 3.7 can be proved for $D$:

**Theorem 5.1.** *For every* $\alpha \in (0, \omega_1)$*, there exists a countable complete metric space* $M_\alpha$ *such that* $D(\mathcal{F}(M_\alpha)) > \alpha$*.*

As a consequence, given the properties enjoyed by $D$, Corollaries 4.1 and 4.2 follow immediately from this theorem.

Of course, since $D(X) \geq \Phi(X)$ for every Banach space $X$, the countably branching diamond graphs can play the role of the $M_\alpha$ spaces. But actually, we can use much simpler sets: let $M_1$ be the graph consisting of two extremities $x_1$ and $y_1$ at distance 2, and of a point $z_1$ at distance 1 from each extremity. As in $D_1$, we put an edge between two vertices if and only if they are at a distance 1 from each other, and the distance $d_1$ on $M_1$ corresponds to be the shortest path metric in a graph. Then we define inductively the metric space $M_\alpha$ for any ordinal $\alpha \in (0, \omega_1)$ in the same way as $D_\alpha$, but starting with $M_1$ instead of $D_1$: if $\alpha = \beta + 1$ is a successor ordinal, $M_\alpha$ is obtained by replacing each edge of $M_1$ by a scaled-down copy of $M_\beta$, actually by $(M_\beta, \frac{d_\beta}{2})$ where $d_\beta$ stands for the distance in $M_\beta$. We still write $x_\alpha$ and $y_\alpha$ for the extremities of $M_\alpha$, and we denote by $z_\alpha$ the point corresponding to $z_1$ in $M_1$. If $\alpha$ is a limit ordinal, $M_\alpha$ is obtained by glueing together the $M_\beta$'s at their extremities.

The idea behind Theorem 5.1 is even simpler than for the $D_\alpha$ spaces:

*Proof of Theorem 5.1.* Assume $\alpha$ is a successor ordinal, and let us consider $S$ an open slice containing $m_{x_\alpha,y_\alpha}$. We can write $m_{x_\alpha,y_\alpha} = \frac{1}{2}m_{x_\alpha,z_\alpha} + \frac{1}{2}m_{z_\alpha,y_\alpha}$ so since the complement of $S$ is convex, $S$ contains either $m_{x_\alpha,z_\alpha}$ or $m_{z_\alpha,y_\alpha}$; say $m_{x_\alpha,z_\alpha}$ without loss of generality. Then

$$\begin{aligned}\|m_{x_\alpha,y_\alpha} - m_{x_\alpha,z_\alpha}\|_{\mathcal{F}(M_\alpha)} &= \left\| \frac{\delta(x_\alpha) - \delta(y_\alpha)}{2} - (\delta(x_\alpha) - \delta(z_\alpha)) \right\| \\ &= \frac{1}{2}\, \|2\delta(z_\alpha) - \delta(x_\alpha) - \delta(y_\alpha)\| \\ &\geq 1,\end{aligned}$$

where the last inequality follows from the existence of a Lipschitz map $f \in \mathrm{Lip}_0(M_\alpha)$ of norm one satisfying $f(x_\alpha) = f(y_\alpha) = 1$ and $f(z_\alpha) = 0$. So $\mathrm{diam}(S \cap B_{\mathcal{F}(M_\alpha)}) \geq 1$.

Notice that in particular, we proved that for every successor ordinal $\alpha$ and for every slice $S$ containing $m_{x_\alpha,y_\alpha}$, there exists $\gamma_S \in \{m_{x_\alpha,z_\alpha}, m_{z_\alpha,y_\alpha}\}$ in $S \cap B_{\mathcal{F}(M_\alpha)}$ such that $\|\gamma_S - m_{x_\alpha,y_\alpha}\| \geq 1$. Now an easy transfinite induction enables us to show that for every $\alpha \in (0, \omega_1)$, we have $m_{x_\alpha,y_\alpha} \in \Delta_\varepsilon^\alpha(B_{\mathcal{F}(M_\alpha)})$ with $\varepsilon = 1$, where the $\Delta_\varepsilon^\alpha(.)$ are the derived sets associated with the dentability index $D$. This finishes the proof. □

Recall that for every Banach space $X$, $\Phi(X) \leq D(X)$. The same inequality is satisfied by the Szlenk index $S_Z$ and the $w^*$-dentability index $D_Z$, and furthermore we know that there exists a function $\Psi\colon (0, \omega_1) \to (0, \omega_1)$ such that for every Banach space $X$ and for any $\alpha \in (0, \omega_1)$, $S_Z(X) \leq \alpha$ implies $D_Z(X) \leq \Psi(\alpha)$ (see [20], Theorem 9). However, for the fragmentability and dentability indices, it is known that the conditions $\Phi(X) < \omega_1$ and $D(X) < \omega_1$ are not even equivalent: the predual $B$ of the James tree space is separable with the PCP so $\Phi(B) < \omega_1$, but does not have the RNP so $D(B) = \infty$. Yet, since these two conditions are equivalent for separable Lipschitz-free spaces, we may ask the following question:

**Question 5.2.** Does there exist a function $\Psi\colon (0, \omega_1) \to (0, \omega_1)$ such that for every separable metric space $M$ and for any $\alpha \in (0, \omega_1)$, $\Phi(\mathcal{F}(M)) \leq \alpha$ implies $D(\mathcal{F}(M)) \leq \Psi(\alpha)$?

Here, since $D(\mathcal{F}(M_\alpha))$ goes increasingly to $\omega_1$ as $\alpha$ tends to $\omega_1$, it makes sense to investigate the behavior of $(\Phi(\mathcal{F}(M_\alpha)))_{\alpha<\omega_1}$. First, we have:

**Lemma 5.3.** *The molecule $m_{x_\omega,y_\omega}$ belongs to $\sigma_1(B_{\mathcal{F}(M_\omega)})$.*

*Proof.* Let $V$ be a $w$-open neighborhood of $m_{x_\omega,y_\omega}$ in $\mathcal{F}(D_\omega)$. For every $k \in \mathbb{N}$, we can write $m_{x_\omega,y_\omega} = \frac{1}{2}m_{x_\omega,z_k} + \frac{1}{2}m_{z_k,y_\omega}$, so using the same arguments as in proof of Proposition 3.1, there exists $k \in \mathbb{N}$ such that $\gamma_V := \frac{1}{2}m_{x_\omega,z_{k+1}} + \frac{1}{2}m_{z_k,y_\omega}$ belongs to $V \cap B_{\mathcal{F}(M_\omega)}$. But

$$\begin{aligned}\|m_{x_\omega,y_\omega} - \gamma_V\|_{\mathcal{F}(M_\omega)} &= \left\| \frac{m_{x_\omega,z_k} - m_{x_\omega,z_{k+1}}}{2} \right\| \\ &= \left\| m_{z_{k+1},z_k} \right\| = 1,\end{aligned}$$

and hence $\mathrm{diam}(V \cap B_{\mathcal{F}(M_\omega)}) \geq 1$, which concludes the proof. □

Actually, the key idea behind this lemma is that $D_1$ embeds isometrically into $M_\omega$ with $t_1$, $b_1$ respectively sent on $x_\omega$, $y_\omega$, and the middle points $(x_1^n)_{n\in\mathbb{N}}$ sent on distinct $z_k$ points. We can generalize this idea: given a limit ordinal $\alpha = \omega.\beta$ where $\beta \in (0, \omega_1)$, it is easy to check that $D_\beta$ embeds isometrically into $M_\alpha = M_{\omega.\beta}$ with $t_\beta$, $b_\beta$ respectively sent on $x_\alpha$, $y_\alpha$, and, if $\beta$ is a successor ordinal, the middle points $(x^n_\beta)_{n\in\mathbb{N}}$ sent on distinct $z_\mu$ points with $\mu < \alpha$. Through this isometric embedding, the molecule $m_{x_\alpha,y_\alpha}$ identifies with $m_{t_\beta,b_\beta} \in \sigma_1^\beta(B_{\mathcal{F}(D_\beta)}) \subsetneqq \sigma_1^\beta(B_{\mathcal{F}(M_{\omega.\beta})})$, and we can deduce:

**Proposition 5.4.** *For every limit ordinal $\alpha = \omega.\beta$ where $\beta \in (0, \omega_1)$, the molecule $m_{x_\alpha,y_\alpha}$ belongs to $\sigma_1^\beta(B_{\mathcal{F}(M_\alpha)})$. In particular, $\Phi(\mathcal{F}(M_\alpha)) > \beta$.*

Hence, $\Phi(\mathcal{F}(M_\alpha))$ also goes increasingly to $\omega_1$ as $\alpha$ tends to $\omega_1$, which does not provide a negative answer to Question 5.2. Thus, this question remains open.

5.2. **About Corollary 4.3.** Corollary 4.3 could also be deduced combining results in [2] and [6]. Indeed, recall Remark 3.8: in [6], Theorem 4.3 states that there exists a uniformly discrete metric space $M$ such that for each Banach space $X$ whose dual contains an isomorphic copy of $\mathcal{F}(M)$, the Szlenk index $S_Z(X)$ of this space is greater than $\omega^2$. In addition, the authors showed that $\mathcal{F}(M)$ is isomorphic to a subspace of a separable dual, so $\mathcal{F}(M)$ has the RNP. Thus, if there exists some compact metric space $K$ such that $\mathcal{F}(M) \simeq \mathcal{F}(K)$, $\mathcal{F}(K)$ has the RNP, which implies that $K$ is a purely 1-unrectifiable compact, and then $\mathcal{F}(K) \cong \mathrm{lip}_0(K)^*$ thanks to Theorem 3.2 in [2]. So the above result implies that the Szlenk index of $\mathrm{lip}_0(K)$ is at least $\omega^2$. But we must have $S_Z(\mathrm{lip}_0(K)) \leq \omega$ since, as mentioned before, $\mathrm{lip}_0(K)$ is AUS renormable. In conclusion, we have:

**Proposition 5.5.** *There exists a uniformly discrete metric space whose free space is not isomorphic to a free space over any compact metric space.*

**Acknowledgments.** This work was partially supported by the French ANR project No. ANR-20-CE40-0006. The author is grateful to Gilles Lancien and Antonín Procházka for their advice, comments and fruitful conversations. The author would also like to thank Ramón J. Aliaga for drawing our attention to some details in the proof of Lemma 3.4.

Université Marie et Louis Pasteur, CNRS, LmB (UMR 6623), F-25000 Besançon, France
*Email address*: estelle.basset@math.cnrs.fr